\documentclass[10pt]{article}
\usepackage{amsmath}
\usepackage{amssymb}
\usepackage{amsthm}
\usepackage{graphicx}
\usepackage{subfigure}
\usepackage{verbatim}
\usepackage{enumerate}
\usepackage{showlabels}

\font\elevenss=cmss11

\font\eightss=cmss8

\font\sixss=cmss8 at 6pt

\newfam\ssfam
\textfont\ssfam=\elevenss \scriptfont\ssfam=\eightss
 \scriptscriptfont\ssfam=\sixss

\def\myffrac#1#2 in #3{\raise 2.6pt\hbox{$#3 #1$}\mkern-1.5mu\raise 0.8pt\hbox{$#3/$}\mkern-1.1mu\lower 1.5pt\hbox{$#3 #2$}}

\newcommand{\Var}{{\text{\rm Var}}}

\theoremstyle{plain}
\newtheorem{thm}{Theorem}[section]
\newtheorem{lem}[thm]{Lemma}
\newtheorem{pr}[thm]{Proposition}
\newtheorem{cor}[thm]{Corollary}

\newtheorem{question}{Question}
\newtheorem*{qth}{Theorem}
\newtheorem*{qlem}{Lemma}
\newtheorem*{qcor}{Corollary}

\theoremstyle{remark}

\newtheorem*{unremark}{Remark}

\def\ee{\epsilon}

\def\vv{{\bf v}}

\def\xx{{\bf x}}

\def\zz{{\bf z}}

\def\aa{{\bf a}}

\def\ww{{\bf w}}
\def\rr{{\bf r}}
\def\XX{{\bf X}}

\def\one{{\bf 1}}

\def\nbd{{\cal N}}
\def\B{{\cal B}}
\def\Cox{\hfill \Box}
\def\C{{\mathbb C}}
\def\Z{{\mathbb Z}}

\def\R{{\mathbb R}}
\def\H{{\mathbb H}}

\def\P{{\mathbb P}}
\def\E{{\mathbb E}}
\def\PP{{\mathcal P}}
\def\F{{\mathcal F}}
\def\nbd{{\mathcal N}}
\def\cL{{\mathcal L}}
\def\Var{{\rm Var}\,}
\def\Bin{{\rm Bin}\,}

\def\Arg{{\rm Arg\,}}

\def\NR{{\rm NR}\,}
\def\disp{\displaystyle}
\def\nullsp{{\cal N}}
\def\ZZZ{{\bf Z}}
\def\cM{{\cal M}}

\begin{document}

\begin{titlepage}
\begin{center}
{\large \bf Multivariate CLT follows from strong Rayleigh property} \\
\end{center}
\vspace{5ex}
\begin{flushright}
Subhroshekhar Ghosh\footnote{Princeton Department of Mathematics, Fine Hall, 
Washington Road, Princeton, NJ 08544, {\tt subhrowork@gmail.com}} \\
Thomas M.\ Liggett\footnote{UCLA Department of Math, 520 Portola Plaza, 
Los Angeles, CA 90095, {\tt tml@math.ucla.edu}} \\
Robin Pemantle\footnote{University of Pennsylvania, 
Department of Mathematics, 209 S. 33rd Street, Philadelphia,
PA 19104 USA, {\tt pemantle@math.upenn.edu}}
\end{flushright}

\vfill

\noindent{\bf ABSTRACT:} 
Let $(X_1 , \ldots , X_d)$ be random variables taking nonnegative
integer values and let $f(z_1, \ldots , z_d)$ be the probability 
generating function.  Suppose that $f$ is real stable; equivalently,
suppose that the polarization of this probability distribution
is strong Rayleigh.  In specific examples, such as occupation counts
of disjoint sets by a determinantal point process, it is 
known~\cite{soshnikov02} that the joint distribution must 
approach a multivariate Gaussian distribution.  We show that 
this conclusion follows already from stability of $f$.
\hfill \\[1ex]

\vfill



\end{titlepage}

\setcounter{equation}{0}
\section{Introduction} \label{sec:outline} 

Let $(X_1 , \ldots , X_d)$ be a random vector on a probability
space $(\Omega , \F , \P)$, whose values lie in a bounded subset
$S$ of the nonnegative orthant in the integer lattice $\Z^d$.
We employ boldface notation for vectors, e.g., $\XX := (X_1 , \ldots , X_d)$
and the monomial product notation $\zz^\rr := z_1^{r_1} \cdots z_d^{r_d}$.
Define the probability generating function $f : \C^d \to \C$ by
$f(\zz) := \sum_{\rr \in S} \P(\XX = \rr) \zz^\rr$.  This is a 
polynomial and is real valued on real inputs.  We say that $f$
is {\em real stable} if $f$ has no zeros in $\H^d$ where $\H$
is the open upper half plane $\{ z : \Im \{ z \} > 0 \}$.  There
is a wealth of knowledge about the behavior of the coefficients
of a real stable polynomial $f$ and, in particular, about the 
resulting probability inequalites when $f$ is a probability 
generating function.  

In the case where $S = \{ 0 , 1 \}^d$ is the $d$-dimensional hypercube, 
the variables $X_j$ are binary valued and real stability of the 
generating function is known as the {\em strong Rayleigh property}.  
Stability is preserved under the following pair of inverse operations:
polarization and aggregation.  If $\{ X_{k,i} : k \leq d, i \leq n_k \}$
is a finite family of nonnegative integer variables with real stable
generating polynomial, the aggregate variables 
$X_k := \sum_{i=1}^{n_k} X_{k,i}$ will also have a real stable
generating polynomial.  This is because it follows from the definition
that stability is preserved by substituting $z_{k,i} := z_k$
for all $k \leq d$ and $i \leq n_k$.  Conversely, if 
$\{ X_k : 1 \leq k \leq d \}$ are random variables whose joint law
has a real stable generating polynomial, $f$, one can define the 
{\em polarization} of $f$ by the substitutions 
$z_k^j = \binom{n_k}{j}^{-1} e_j (z_{k,1} , \ldots , z_{k,n_k})$ 
where $n_k$ are upper bounds for the values of $X_k$ and $e_j$ is the 
elementary symmetric function of degree $j$ on $n_k$ variables.
Stability of $f$ implies stability of the polarization $\PP f$
of $f$, hence the strong Rayleigh property for a collection of binary
variables $\{ X_{k,i} \}$ with probability generating function $\PP f$.
This result may be found in~\cite[Theorem~4.7]{borcea-branden-liggett}
and follows from the Grace-Walsh-Szeg{\"o} coincidence Theorem.

When $d=1$, the theory of distributions with real stable generating
polynomials is completely understood.  If $f(z)$ has no zeros in $\H$
then, by invariance under conjugation, it has only real zeros.  By 
nonnegativity of the coefficients, all zeros of $f$ lie in the
negative half line.  Writing $f(z) = C \prod_{j=1}^d (z + a_j)$,
dividing the $j^{th}$ term by $1+a_j$ and using the fact that
$f(1) = 1$, we have the representation
$$f(z) = \prod_{j=1}^d (1 - p_j + p_j z)$$
where $p_j := 1/(1+a_j)$ are numbers in the unit interval.  In other
words, a variable $X$ with generating function $f$ is distributed as 
a sum of indepdendent Bernoulli variables with means $p_1, \ldots , p_d$.
This implies a self-normalized central limit theorem: 
\begin{equation} \label{eq:one var}
\frac{X - \E X}{(\Var X)^{1/2}} \to N(0,1) 
   \;\; \mbox{ as } \;\; \Var (X) \to \infty \, .
\end{equation}

When $d \geq 2$, real stable polynomials do not in general factor; 
therefore this argument cannot be applied to establish a multivariate 
central limit theorem.  Nevertheless, such a result is known to
hold in one of the most important applications of strong Rayleigh
distributions, namely occupations of determinantal processes.
To elaborate, let $N : \Omega \times \B \to \Z^+$ be a random
counting measure on the Borel subsets of $\R^n$ and suppose that
$N(\cdot ,\cdot)$ is Hermitian determinantal.  This means that the $k$-fold 
joint intensities exist for all $k$ and the joint intensity at 
$(\xx^{(k)} , \ldots , \xx^{(k)})$ is the determinant of the matrix 
$\disp \left [ K(\xx^{(i)} , \xx^{(j)})_{1 \leq i,j \leq k} \right ]$
where $K$ is a Hermitian kernel.  Let $B_1 , \ldots , B_d$ be disjoint
Borel subsets of $\R^n$.  It is known that the counts
$X_j := N(B_j)$ for $1 \leq j \leq d$ have real stable generating
polynomial provided that each $N(B_j)$ is bounded; this follows,
for example from the determinantal form of the generating function
given in~\cite[Theorem~2]{soshnikov00} together with 
with~\cite[Proposition~3.2]{borcea-branden-liggett}.
Boundedness is not an overly restrictive assumption because under
finiteness of the means $\E N(B_j)$ one can always approximate 
the occupations $N(B_j)$ in total variation by bounded variables.
Soshnikov~\cite[p.~174]{soshnikov02} proved a normal limit theorem
for linear combinations $\sum_{j=1}^d \alpha_j X_j$, which is equivalent
to a multivariate CLT.  This generalized an earlier result for
several specific determinantal kernels arising in random 
spectra~\cite{soshnikov00}.

Determinantal measures are in some sense a very small set of measures.  
For example, determinantal measures supported on a set of cardinality $d$
are parametrized by $d \times d$ Hermitian matrices, and therefore
occupy a $d^2$-dimensional set in the $(2^d - 1)$-dimensional space
of probability laws on $\{ 0 , 1 \}^d$.  The set of strong Rayleigh
measures, by contrast, has full dimension, being constrained by 
inequalities rather than identities.  Because of the relative robustness
of the strong Rayleigh property, it seems useful to discover whether 
properties of determinantal measures, such as multivariate Gaussian 
behavior, follow already from stability.  

Our main results, Theorem~\ref{th:first} in the bivariate case and 
Theorem~\ref{th:first}$'$ in the multivariate case, show this to be 
the case.  The subsequent sections discuss extensions and some 
theoretical questions about the class of real stable distributions
which are raised by the arguments of the paper and partially answered.

\setcounter{equation}{0}
\section{Main result}

Our first result in this direction is a bivariate CLT valid
when the variance grows faster than the $2/3$ power of the
maximum value.  Because real stable variables are known to 
be negatively correlated, the covariances are denoted by
negative quantities. 

\begin{thm} \label{th:first}
Let $\{ (X_n , Y_n) : n \geq 1 \}$ be a sequence of random integer
pairs each of whose bivariate generating polynomials $f_n(x,y)$ 
is real stable and has degree at most $M_n$ in each variable.  Let
$$A_n = \left [ \begin{array}{cc} \alpha_n & -\beta_n \\ -\beta_n & \gamma_n 
   \end{array} \right ]$$
denote the covariance matrix of $(X_n , Y_n)$.  Suppose there is
a sequence $s_n \to \infty$ and a fixed matrix 
$\disp A = \left [ \begin{array}{cc} \alpha & -\beta \\ 
-\beta & \gamma \end{array} \right ]$
such that $s_n^{-2} A_n \to A$ and $s_n^{-1} M_n^{1/3} \to 0$.  Then
\begin{equation} \label{eq:statement}
\frac{(X_n , Y_n) - (\E X_n , \E Y_n)}{s_n} \to N(0 , A)
\end{equation}
in distribution as $n \to \infty$.
\end{thm}

An outline of the proof is as follows.  Let $a$ and $b$ be positive 
integers.  From the definition of stability it may be shown that 
the generating polynomial for $a X_n + b Y_n$ has no zeros near~1
(this is Lemma~\ref{lem:aX+bY} beow).  This implies a Gaussian 
approximation for $a X_n + b Y_n$ (Lemma~\ref{lem:sufficient CLT} 
below).  Tightness and continuity could be used to extend this to 
positive real $(a,b)$, however the usual Cram{\'e}r-Wold argument 
requires this for all real $(a,b)$ regardless of sign.  Instead, the
argument is finished instead by invoking an improved Cram{\'e}r-Wold 
result (Lemma~\ref{lem:cramer-wold} and Corollary~\ref{cor:cramer-wold}).

\begin{lem} \label{lem:aX+bY}
Whenever $(X,Y)$ is stable and $b \geq a$ are positive integers,
the probability generating function for $aX+bY$ has no zeros
in the open disk of radius $\delta$ about~1, where
$\delta := \sin (\pi/b)$.
\end{lem}

\noindent{\sc Proof:} If $f(x,y)$ is the pgf for $(X,Y)$ then
the pgf for $aX+bY$ is $f(z^a,z^b)$.  Stability of $f$ implies
that $f(z^a,z^b)$ has no zeros whose argument $z$ lies in the
open interval $(0 , \pi/b)$.  Invariance under conjugation 
and the fact that a probability generating function can never 
have positive real zeros implies that $f(z^a,z^b)$ is in fact
zero-free on the sector $\{ z : |\Arg (z)| < \pi / b \}$.  
The nearest point to~1 in this sector is at distance $\delta$.  
$\Cox$

\begin{lem}[\protect{\cite[Corollary~4.3]{cramer-wold}}]  
\label{lem:cramer-wold}
Let $\cL$ be an infinite family of $(d-1)$-dimensional subspaces 
of $\R^d$.  Let $\pi_L$ denote projection of measures onto $L$,
in other words $\pi_L \mu := \mu \circ \pi_L^{-1}$.  Let $\mu$ be 
a probability measure on $\R^d$ with finite moment generating 
function in a neighborhood of the origin and let $\nu$ be any 
probability measure on $\R^d$.  Suppose that the projections 
$\pi_L \mu$ and $\pi_L \nu$ coincide for every $L \in \cL$.  
Then $\mu = \nu$.
$\Cox$
\end{lem}

\begin{cor} \label{cor:cramer-wold}
Let $\mu$ be a centered Gaussian law on $\R^d$ and let $\cL$ be 
an infinite family of $(d-1)$-dimensional subspaces of $\R^d$.  
Suppose $\{ \mu_n \}$ is a sequence of probability measures
on $\R^d$ such that for each $L \in \cL$, the projections
$\pi_L \mu_n$ converge in the weak topology as $n \to \infty$ 
to $\pi_L \mu$.  Then $\mu_n \to \mu$.
\end{cor}

\noindent{\sc Proof:} 
Convergence of $\pi_L \mu_n$ for more than one hyperplane $L$
implies tightness of the family $\{ \mu_n \}$.  Therefore,
any subsequence of $\{ \mu_n \}$ has a convergent sub-subsequence;
denote its limit by $\nu$.  It suffices to show that $\nu = \mu$.
Each $\pi_L$ is continuous, therefore $\pi_L \nu = \lim_{n \to \infty} 
\pi_L \mu_n = \pi_L \mu$.  Noting that $\mu$ has moment generating 
function defined everywhere, the conclusion now follows from 
Lemma~\ref{lem:cramer-wold}. 
$\Cox$

\begin{lem}[\protect{\cite[Theorem~2.1]{LPRS2015}}] 
\label{lem:sufficient CLT}
Let $f$ be the generating polynomial for a probability law $Q$ on
the nonnegative integers.  Let $N$ denote the degree of $f$.
Let $m$ and $\sigma^2$ respectively
denote the mean and variance of $Q$ and let $F$ denote the 
self-normalized cumulative distribution function defined by 
$$F(x) := \sum_{k \leq m + x\sigma} Q(k) \, .$$
Let $\nbd (x) := (2\pi)^{-1/2} 
\int_{-\infty}^x e^{-t^2/2} \, dt$ denote the standard normal CDF. 
Given $\delta > 0$, there exists a constant $C_\delta$ depending 
only on $\delta$ such that if $f$ has no roots in the ball
$\{ z : |z - 1| < \delta \}$ then 
$$\sup_{x \in \R} |F(x) - G(x)| \leq 
   C_\delta \frac{N^{1/3}}{\sigma} \, .$$
\end{lem}

\noindent{\sc Proof:} The result as stated in~\cite[Theorem~2.1]{LPRS2015},
in the special case $z_0 = 1$ has the upper bound $B_1 N/\sigma^3
+ B_2 N^{1/3}/\sigma$ with $B_1$ and $B_2$ depending on $\delta$.
Because $|F - \nbd|$ is never more than~1, we may assume that 
$N^{1/3}/\sigma \leq B_2^{-1}$, whence $B_1 N/\sigma^3 \leq 
(B_1 / B_2^2) N^{1/3}/\sigma$.  Setting $C = C_\delta = 
B_2 + B_1 / B_2^2$ recovers the result in our form.  The
result as stated holds for $N > N_0 (\delta)$, but with
$C N_0$ in place of $C$ it holds for all $N$.  
$\Cox$

\noindent{\sc Proof of Theorem}~\ref{th:first}:  
We will apply Corollary~\ref{cor:cramer-wold} with $\mu = N(0,A)$
and $\cL$ equal to the set of lines through the origin with 
positive rational slope.  Given $L \in \cL$, let $(a,b)$ be
a positive integer pair in $L$.  Then $\pi_L (X,Y) = 
(a X + b Y) / \sqrt{a^2+b^2}$ and $\pi_L \mu = N(0,V)$ where 
$$V := V (a , b) := \frac{\alpha \, a^2 - 2\beta \, a b + \gamma \, b^2} 
   {a^2 + b^2} \, .$$  
According to Corollary~\ref{cor:cramer-wold}, the theorem will follow 
if we can show that 
\begin{equation} \label{eq:need}
\frac{a}{\sqrt{a^2 + b^2}} \frac{X_n - \E X_n}{s_n} 
   + \frac{b}{\sqrt{a^2 + b^2}} \frac{Y_n - \E Y_n}{s_n} 
   \to N (0 , V(a,b))
\end{equation}
weakly for fixed positive integers $a$ and $b$ as $n \to \infty$.
We proceed to show this.

First, if $V(a,b) = 0$, we observe that the left-hand side 
of~\eqref{eq:need} has mean zero and variance 
$$\frac{\alpha_n \, a^2 - 2\beta_n \, a b + \gamma_n \, b^2} 
   {(a^2 + b^2) s_n^2} = o(1)$$  
by the assumption that $A_n / s_n^2 \to A$.  Weak convergence to
$\delta_0$, which is the right-hand-side of~\eqref{eq:need}, 
follows from Chebyshev's inequality.  

Assume now that $V \neq 0$.  
By Lemma~\ref{lem:aX+bY}, for all $n$, the generating polynomial $g_n$ 
for $aX_n + bY_n$ has no zeros within distance $\delta := \sin (\pi / b)$ 
of~1.  Apply Lemma~\ref{lem:sufficient CLT} to the generating 
polynomial $g_n$ with $N = (a+b) M_n$.  
In the notation of Lemma~\ref{lem:sufficient CLT},
\begin{eqnarray*}
m & = & a \E X_n + b \E Y_n \, ; \\
\sigma^2 & = & a^2 \alpha_n - 2 a b \beta_n + b^2 \gamma_n \, .
\end{eqnarray*}
The assumption $s_n^{-2} A_n \to A$ implies that
\begin{equation} \label{eq:sigma}
\sigma^2 / s_n^2 \to (a^2 + b^2) V \, .
\end{equation}
The conclusion of the lemma is that $[a(X_n - \E X_n) + b(Y_n - \E Y_n)] /
\sigma$ differs from a standard normal by at most
$$ C_\delta \frac{N^{1/3}}{\sigma} = (1 + o(1)) C_\delta 
   \frac{(a+b)^{1/3}}{\sqrt{a^2 + b^2} \sqrt{V}} \frac{M_n^{1/3}}{s_n} \, .$$
By hypothesis $M_n^{1/3} / s_n \to 0$.  Thus 
$[a(X_n - \E X_n) + b(Y_n - \E Y_n)] / \sigma \to N(0,1)$;
multiplying through by $\sqrt{V}$ and plugging in~\eqref{eq:sigma}
gives~\eqref{eq:need}.
$\Cox$

\begin{question}
Can the hypothesis $s_n^{-1} M_n^{1/3} \to 0$ be weakened, 
preferably to $s_n \to \infty$?
\end{question}

\setcounter{equation}{0}
\section{Extensions}

\subsubsection*{Higher dimensions}

More or less the same argument works to prove a multivariate CLT 
for real stable distributions in $d$ variables.  It requires 
only small generalizations of two lemmas, the first of which
is immediate. 
\begin{qlem}[${\bf \protect{\ref{lem:aX+bY}}'}$]
If $(\XX)$ is stable and $\aa$ is a positive integer vector, 
then the probability generating function for $\aa \cdot \xx$
has no zeros in the open disk of radius $\delta$ about~1, where
$\delta := \sin (\pi/\max_j a_j)$.
$\Cox$
\end{qlem}

\begin{qcor}[${\bf \protect{\ref{cor:cramer-wold}}'}$]
Let $\mu$ be a centered Gaussian law on $\R^d$.  If 
$\{ \mu_n \}$ is a sequence of probability measures
such that for all positive rational lines $L$, the projections
$\pi_L \mu_n$ converge to $\pi_L \mu$, then $\mu_n \to \mu$.
\end{qcor}

\noindent{\sc Proof:}
We prove by induction on $m$ that $\pi_L \mu_n \to \pi_L \mu$
for all $m$-dimensional subspaces containing a positive rational 
point.  It is true by hypothesis when $m=1$.  Assume for induction
that it is true for dimensions smaller than $m$. 
Fix any $m$-dimensional subspace $L$ containing a
positive rational point and apply Corollary~\ref{cor:cramer-wold} 
with $\pi_L \mu$ in place of $\mu$ and the infinite family of 
subspaces $L' \subseteq L$ having a basis of $m-1$ positive rational 
vectors in place of $\cL$.  By the induction hypothesis, each 
$\pi_{L'} \mu_n \to \pi_{L'} \mu$, so by Corollary~\ref{cor:cramer-wold},
$\pi_L \mu_n \to \pi_L \mu$, completing the induction.  Once $m=d$,
the corollary is proved.
$\Cox$

These two results imply the extension of Theorem~\ref{th:first}
to $d$ variables:
\begin{qth}[${\bf \protect{\ref{th:first}}}'$]
Let $\{ \XX^{(n)} \}$ be a sequence of random vectors with real stable
generating polynomials, degree at most $M_n$ in each variable, and
covariance matrices $A_n$.  Suppose $s_n \to \infty$ with 
$s_n^{-2} A_n \to A$ and $s_n^{-1} M_n^{1/3} \to 0$.  Then
$(\XX - \E \XX) / s_n \to N(0,A)$ in distribution as $n \to \infty$.
$\Cox$
\end{qth}

\subsubsection*{Singularity of $A$}

When $A$ is singular, say $\langle a,b \rangle A = 0$, the conclusion of 
Theorem~\ref{th:first}, namely a bivariate Gaussian limit, implies 
only that $(aX+bY)/s_n \to 0$, not that $aX+bY$ has a normal limit.
This can be improved to the following result.

\begin{thm} \label{th:singular}
In the notation of Theorem~\ref{th:first}$'$, suppose $A$ is 
singular and let $\nullsp$ denote the nullspace of $A$.  Let
$\one_G$ denote the vector whose $j$ component is~1 if $j \in G$
and~0 otherwise.  The space $\nullsp$ is spanned by a collection
$\{ 1_G : G \in \cM \}$ where $\cM$ is a collection of disjoint sets.
The quantities $Z_G^{(n)} := \one_G \cdot \XX^{(n)}$ all have 
normal limits, provided the variances $\sigma_G^{(n)} := 
\Var (Z_G^{(n)})^{1/2}$ go to infinity; assuming this,  
$(\sigma_G^{(n)})^{-1} (Z_G^{(n)} - \E Z_G^{(n)}) \to N(0,1)$.
\end{thm}

\begin{unremark}
This gives a CLT for a collection of linear functionals spanning
the null space of $A$.  More generally, one might want a CLT for every
element of the null space.  If the null space has dimension $r$ then
one may construct $\{ Z_1 , \ldots , Z_r \}$ as above.  The vectors
$\{ \ZZZ^{(n)} \}$ are real stable with covariance matrices $A_n'$
for which $s_n^{-2} A_n' \to 0$.  If it is possible to find $s_n'$ 
for which $(s_n')^{-2} A_n' \to A'$ then one obtains a finer 
multivariate CLT.  The covariance matrices $A_n'$ may or may not 
have a rescaled limit.
\end{unremark}

\begin{lem} \label{lem:nonsingular}
If $(X_1 , \ldots , X_r)$ is a random integer vector whose 
$r$-variate generating function is real stable, then its 
covariance matrix has nonnegative row and column sums.  
\end{lem}

\noindent{\sc Proof:}  The row sums of the covariance matrix are
the values $\E (X_i - \mu_i) \sum_j (X_j - \mu_j)$.  The argument 
may be reduced to
the case $r=2$ by considering the pair $(X_i , Y_i)$ where 
$Y_i := \sum_{i \neq j} X_j$.  Without loss of generality, we
therefore assume $r=2$ and denote the pair $(X_1, X_2)$ by $(X,Y)$.
We first claim that for all $k$,
\begin{equation} \label{eq:SCP}
\E (Y | X = k) \leq \E (Y | X = k-1) \leq \E (Y | X = k) + 1 \, .
\end{equation}
This follows from the strong Rayleigh property for the polarization
$(X_1 , \ldots , X_m , Y_1 , \ldots , Y_m)$ of $(X,Y)$.  Indeed,
suppose the polarization is coupled to $(X,Y)$ so that for each $k , \ell$,
the conditional law of $(X_1 , \ldots , X_m , Y_1 , \ldots , Y_m)$ 
given $X=k, Y=\ell$ is the product $\nu_{m,k} \times \nu_{m,\ell}$
where $\nu_{m,j}$ is uniform on sequences of zeros and ones of
length $m$ summing to $j$.  Then 
$$\E (Y | X = k) = \E  \left ( \left. \sum_{j=1}^m Y_j \right |
   X_1 = \cdots = X_k = 1 , X_{k+1} = \cdots = X_m = 0 \right )$$
and the claim follows from the stochastic covering property
for strong Rayleigh measures~\cite[Proposition~2.2]{PP-rayleigh}.
In fact it is only the right-hand inequality of~\eqref{eq:SCP} that we 
need.  Adding $X$ gives $\E (X+Y | X = k-1) \leq \E (X+Y | X = k)$.  
Thus $\E (X+Y | X)$ is a monotone increasing function of $X$.  This
immediately implies nonnegative correlation of the bounded variables
$X$ and $X+Y$, which is the conclusion of the lemma.
$\Cox$

The next lemma is stated generally though it is used for one specific
purpose, namely for the covarinace matrix of a collection of 
random integers with real stable generating function.

\begin{lem} \label{lem:ones}
Let $M$ be any symmetric matrix with nonnegative diagonal entries, 
nonpositive off-diagonal entries and nonnegative row sums.
Then the index set $[m]$ may be partitioned into disjoint sets 
$T$ and $\{ S_\alpha \}$ such that $M_{i,j} = 0$ when $i$ and $j$
are in different sets of the partition.  This can be done in such
a way that the restriction $M|_T$ is nonsingular, while the 
restrictions $M|_{S_i}$ have one-dimensional null spaces containing
the vectors with all entries equal.
\end{lem}

\noindent{\sc Proof:}
Recall that $\nullsp$ denotes the null space of $M$.
Choose any nonzero vector $\vv \in \nullsp$ with minimal support set $S$, 
meaning that no vector whose support is a proper subset of $S$ is in 
the null space of $M$.  Suppose $\vv$ has coordinates of mixed sign.  
Let $E$ be the set of indices of positive coordinates and $F$ the 
set indices of negative coordinates.  Let $M'$ be the $2 \times 2$
matrix indexed by the set $\{ E , F \}$ whose $(G,G')$-element is
$\sum_{i \in G , j \in G'} M_{ij}$.
This matrix also has nonnegative diagonal entries (follows from 
nonnegativity of row sums and nonpositivity of off-diagonal elements), 
nonpositive off-diagonal entries (obvious) and nonnegaive row sums.
It has a vector of mixed signs in its null space, namely 
$(\sum_{j \in E} v_j , \sum_{j \in F} v_j)$, hence must be the 
$2 \times 2$ zero matrix.  This means that the $\vv_E$ and $\vv_F$ 
are each separately in the null space (where $\vv_G$ denotes the 
vector whose $j^{th}$ coordinate is $v_j \one_G (j)$).  
This contradicts the minimality of the support of $\vv$. 
We conclude that all elements of the null space with minimal 
support have coordinates all of one sign.

Still assuming $\vv$ to have minimal support set $S \subseteq \nullsp$,
consider the sub-collection $\{ X_j : j \in S \}$, which inherits the
properties in the hypotheses.  Its covariance matrix $M'$ is the submatrix 
of $M$ indexed by $S$.  Assume for contradiction that the coordinates of 
$\vv$ are not equal.  Let $\ww$ be the all ones vector of the same length 
as $\vv$.  Scale $\vv$ so that its minimum coordinate is equal to~1.  
If $v_i = 1$ then
$$0 = (M' \vv)_i \geq (M' \ww)_i \geq 0 \, ,$$
the last inequality following from nonnegativity of the row sums.
It follows that $M_{ij} = 0$ for all $i,j$ such that $v_i = 1 < v_j$.
Thus $S' := \{ i : v_i = 1 \}$ is a proper subset of $S$ whose
indicator vector is in the null space of $M$.  By contradiction, 
$\vv = \ww$ as desired.  

Finally, if $\ww_S$ and $\ww_T$ are vectors of ones and zeros with 
support sets $S$ and $T$ respectively and these are not disjoint, 
then $\ww_S - \ww_T \in \nullsp$ and is of mixed sign, a contradiction.  
This finishes the proof.
$\Cox$ 

\noindent{\sc Proof of Theorem}~\ref{th:singular}:
The conclusions of Lemma~\ref{lem:nonsingular} pass to the limit:
the limiting covariance matrix $A$ has nonnegative row sums as well
as being symmetric with nonnegative diagonal entries and nonpositive
off-diagonal entries.  The conclusions of Lemma~\ref{lem:ones} then
follow as well.  Fix $\vv$ such that $\vv A = 0$.  It follows from
Lemma~\ref{lem:ones} that $\ww A = 0$ as well.  The random 
variables $Z_n = \ww \cdot \XX^{(n)}$ are univariate real stable, 
hence subject to the real stable CLT~\eqref{eq:one var}.  In particular,
$\sigma_n^{-1} (Z^{(n)} - \E Z^{(n)}) \to N(0,1)$ 
weakly whenever $\sigma_n := \Var (Z_n)^{1/2} \to \infty$. 
$\Cox$

\subsubsection*{Quantitative version}

Suppose $f_n$ is a sequence of bivariate real stable generating 
functions and that $M_n / s_n^3$ goes to zero, where $M_n$ is
the maximum degree of $f_n$ in either variable and $s_n^2$ is
the maximum variance of either variable.  Let $Q_n$ denote the
probability law represented by $f_n$ and let $A_n$ denote the
covariance matrix for this law.  Suppose that the $Q_n$, centered 
and divided by $s_n$, stays at least $\ee$ away from the bivariate 
Gaussian with mean zero and covariance $s_n^{-2} A_n$.  Taking a
subsequence $\{ n_k \}$, there is a matrix $A$ such that
$s_n^{-2} A_n \to A$, contradicting Theorem~\ref{th:first}.
We conclude that there is a quantitative version of this result:
namely a function $g$ going to zero at zero such that
\begin{equation} \label{eq:quant}
||Q - N(\vv , \Sigma)|| < g(M^{1/3} / ||\Sigma||^{1/2})
\end{equation}
whenever $Q$ is a bivariate real stable law with mean $\vv$, covariance 
$\Sigma$ and maximum $M$.

\begin{question}
What is the best possible function $g$ in~\eqref{eq:quant}?
\end{question}

Lemmas~\ref{lem:aX+bY} and~\ref{lem:sufficient CLT} are quantitative
and sharp.  Therefore, establishing~\eqref{eq:quant} without giving 
up too much in the choice of function $g$ would rest on a quantitative
version of Corollary~\ref{cor:cramer-wold}.  Inverting the 
characteristic function is inherently quantitative, however
the use of uniform continuity so as to use only values on a 
finite mesh is messy.  Furthermore, while Lemma~\ref{lem:aX+bY}
is sharp, its use is certainly not: for example, if $f(z)$ generates
a distribution within $\ee$ of normal, then so does $f(z^k)$,
even though the nearest zero to~1 becomes nearer by a factor of $k$.

Non-uniformity of the estimates as the denominator of the rational
slope increases is an annoying artifact of the proof and points to
the need to replace Lemma~\ref{lem:aX+bY} with something uniform
over sets of directions.  One possibility is to replace the exact
combination $aX+bY$ with $a,b \in \Z^+$ by a probabilistic 
approximation.  One somewhat crude approximation is to let 
$Z := \Bin (X,a) + \Bin (Y,b)$ be the sum of binomial distributions,
conditionally independent given $(X,Y)$.  This has generating
polynomial $g(z) = f(1-a+az,1-b+bz)$ if $f(x,y)$ is the generating
polynonmial for $(X,Y)$.  When $f$ is stable, so is $g$, thereby 
achieving uniformity in direction.  Conditioned on $(X,Y)$, the
difference $Z - aX - bY$ is normal with variance $a(1-a)X + b(1-b)Y$,
which has order $M$.  The size parameter $M$ cannot be less than
a constant times $s^2$, where $s^2$ is the norm of the covariance
matrix, but in the regime where $M = O(s^2)$, the added noise does
not swamp the signal and near normality of $Z$ implies near normality
of the true $aX+bY$.  This works equally well in any dimension.

To extend beyond the regime where $M$ and $s^2$ are comparable, 
we would need to find a random variable $Z$ with real stable law
that approximates $aX+bY$ to within a smaller error than $M^{1/2}$.
This motivates a one-dimensional version of this problem, which
we now discuss.

\setcounter{equation}{0}
\section{Approximate multiplication}

Let $X$ be a positive integer random variable with stable 
generating polynomial $f$.  We use ``stable multiplication by $a$''
to denote the construction of a random variable $Z$ with 
$|Z - aX| = O(1)$.  

\begin{pr}[stable division by~2] \label{pr:half}
Conditional on $X$, if $X$ is even let $Z = X/2$, while if $X$ is odd,
flip a fair coin to decide whether $Z = \lfloor X/2 \rfloor$ or
$Z = \lceil X/2 \rceil$.  Then $Z$ stably multiplies $X$ by $1/2$.
\end{pr}

\begin{pr}[stable division by~$k$] \label{pr:divide}
For any $k \geq 2$, $\lfloor X/k \rfloor$ stably multiplies $X$ by $1/k$.
\end{pr}

The engine for proving both of these is the following result concerning
interlacing roots.
Let $\NR$ be the collection of polynomials all of whose roots are 
simple and strictly negative.  If $f$ is a polynomial of degree $n$ 
and $k\geq 1$, write 
\begin{equation}\label{poly}
f(x)=\sum_{i=0}^{k-1}x^ig_i(x^k),
\end{equation}
where $g_i$ is a polynomial of degree $\lfloor\frac{n-i}k\rfloor$.

\begin{thm} \label{th:interlace}
If $f\in \NR$ has degree $n$, the corresponding polynomials 
$g_i$ are in $\NR$ as well. Furthermore, their roots are interlaced 
in the sense that if the collection of all $n-k+1$ roots $s_j$ of the 
$g_i$'s are placed in increasing order,
$$s_{n-k}<\cdots<s_4<s_3<s_2<s_1<s_0<0,$$
then the roots of $g_i$ are $s_i, s_{i+k}, s_{i+2k},\dots.$
\end{thm}

\noindent{\sc Proof:}
The proof is by induction on the degree $n$ of $f$. Let $r_1,\dots,r_n$ 
be the negatives of the roots of $f$, and let $e_j=e_j(r_1,\dots,r_n)$
be the elementary symmetric functions:
$$e_0=1,\quad e_1=\sum_ir_i,\quad e_2=\sum_{i<j}r_ir_j,\dots.$$
Assuming without loss of generality that $f$ is monic,
\begin{equation}\label{product}f(x)=\prod_{i=1}^n(x+r_i)=\sum_{j=0}^nx^je_{n-j}.\end{equation}
Then
$$g_i(y)=\sum_{j=0}^{\lfloor\frac{n-i}k\rfloor}y^je_{n-kj-i}.$$

For the base step of the induction, take $n<2k$, so that the $g_i$'s are linear or constant. In fact, $g_i(y)=e_{n-i}$  if $i>n-k$ and 
$g_i(y)=e_{n-i}+ye_{n-k-i}$ if $i\leq n-k$. In the latter case, the root is $-e_{n-i}/e_{n-k-i}$,
so the interlacement property is a consequence of the log concavity of the sequence $e_m$:
$$\frac{e_{m+1}}{e_m}\downarrow.$$
This statement is a consequence of Newton's inequalities; see~\cite{HLP} 
and~\cite{rosset1989}.

Now assume the result for a given $n$, let $f$ be as in (\ref{product}),  consider the polynomial of degree $n+1$
$$F(x)=(x+r)f(x),\quad r>0,$$
and its decomposition
$$F(x)=\sum_{i=0}^{k-1}x^iG_i(x^k).$$
If $e_j'=e'_j(r_1,\dots,r_n,r)$ are the elementary symmetric functions corresponding to the longer sequence,
$e_j'=e_j+re_{j-1}$, so
\begin{equation}\label{Gg}\begin{gathered}G_i(y)=\sum_{j=0}^{\lfloor\frac{n+1-i}k\rfloor} y^je'_{n+1-kj-i}=\sum_{j=0}^{\lfloor\frac{n+1-i}k\rfloor} y^j[e_{n+1-kj-i}+re_{n-kj-i}]\\
=rg_i(y)+\begin{cases}yg_{k-1}(y)&\text{ if }i=0;\\g_{i-1}(y)&\text{ if }i\geq 1.\end{cases}\end{gathered}\end{equation}
Now we use this to determine the sign of $G_i(s_j)$. The signs of $g_i$ alternate between intervals separated by the roots of $g_i$, since all roots are simple.
Also, $g_i(0)>0$ for each $i$. 

We describe the argument in the following array, in case $k=3$:
$$\left(\begin{matrix}&\cdots&s_{6}&s_{5}&s_{4}&s_{3}&s_{2}&s_{1}&s_0&0\\\\
g_0&\cdots&0&+&+&0&-&-&0&+\\
g_1&\cdots&+&+&0&-&-&0&+&+\\
g_2&\cdots&+&0&-&-&0&+&+&+\\\\
G_0&\cdots&-&+&+&+&-&-&-&+\\
G_1&\cdots&+&+&+&-&-&-&+&+\\
G_2&\cdots&+&+&-&-&-&+&+&+\end{matrix}\right).$$
Note that each row is periodic of period~6, and each row within the 
two groups is obtained from the previous row via a shift. 
Here are some examples of the computation for the bottom rows:
$$G_0(s_2)=rg_0(s_2)+s_2g_2(s_2)=rg_0(s_2)<0, G_2(s_3)=rg_2(s_3)+g_1(s_3)<0.$$

More generally note that the induction hypothesis implies that
\begin{equation}\label{gi}g_i(s_j)\begin{cases}<0&\text{ if }\frac{j-i}k\in \cup_{m=0}^{\infty}(2m,2m+1);\\=0&\text{ if }\frac{j-i}k\in\{0,1,2,\dots\};\\
>0&\text{ if }\frac{j-i}k\in (-\infty,0)\cup\cup_{m=0}^{\infty}(2m+1,2m+2).\end{cases}\end{equation}
We would like to show that
\begin{equation}\label{Gi}G_i(s_j)\begin{cases}<0&\text{ if }\frac{j-i}k\in \cup_{m=0}^{\infty}[2m,2m+1);\\
>0&\text{ if }\frac{j-i}k\in (-\infty,0)\cup\cup_{m=0}^{\infty}[2m+1,2m+2).\end{cases}\end{equation}

There are several cases to consider. First take $i=0, 2mk\leq j<(2m+1)k$ for some $m\geq 0$. Then by (\ref{Gg}),
$$G_0(s_j)=rg_0(s_j)+s_jg_{k-1}(s_j).$$
By (\ref{gi}), $g_0(s_j)=0$ if $j=2mk$ and is $<0$ otherwise, while $g_{k-1}(s_j)=0$ if $j=(2m+1)k-1$ and $>0$ otherwise. Since $r>0$
and $s_j<0$, $G_0(s_j)<0$ as required.  The next case is $i=0, (2m+1)k\leq j<(2m+2)k$ for some $m\geq 0$. Now
$g_0(s_j)=0$ if $j=(2m+1)k$ and $>0$ otherwise, while $g_{k-1}(s_j)=0$ if $j=(2m+2)k-1$ and $<0$ otherwise, so $G_0(s_j)>0$.

Next take $i\geq 1$ and  $2mk\leq j-i<(2m+1)k$ for some $m\geq 0$. Now
$$G_i(s_j)=rg_i(s_j)+g_{i-1}(s_j),$$
$g_i(s_j)=0$ if $j-i=2mk$ and is $<0$ otherwise, and $g_{i-1}(s_j)=0$ if $j-i=2(m+1)k-1$ and is $<0$ otherwise, so $G_i(s_j)<0$.
If, on the other hand, $i\geq 1$ and $j<i$ or $(2m+1)k\leq j-i<(2m+2)k$ for some $m\geq 0$,
$g_i(s_j)=0$ if $j-i=(2m+1)k$ and is $>0$ otherwise, and $g_{i-1}(s_j)=0$ if $j-i=(2m+2)k-1$ and is $>0$ otherwise, so $G_i(s_j)>0$.

From (\ref{Gi}) we see that $G_i$ has a root in each interval of the form 
\begin{equation}\label{interval}(s_{mk+i},s_{mk+i-1})\end{equation} 
for $0\leq m\leq\frac{n-k-i}k$. (By convention, we set $s_{-1}=0$.) This shows that $G_i$ has at least $\lfloor\frac{n-i}k\rfloor$ negative roots.
The degree of $G_i$ is $\lfloor\frac{n+1-i}k\rfloor.$ We see that all roots of $G_i$ are negative, except possibly in case
$\lfloor\frac{n-i+1}k\rfloor=\lfloor\frac{n-i}k\rfloor+1$. In this case, the extra root is recovered by noting that, with $m=\lfloor\frac{n-i}k\rfloor$,
$$(-1)^mG_i(s_{(m-1)k+i})>0\quad\text{and}\quad (-1)^mG_i(s)<0\text{ for large negative }s.$$
Therefore, $G_i$ has the correct number of negative roots. The interlacement property follows from the form of the intervals in (\ref{interval}):
$$t_{n-k+1}<s_{n-k}<t_{n-k}\cdots<s_2<t_2<s_1<t_1<s_0<t_0<0.$$
where the roots of $G_i$ are $t_i,t_{i+k},t_{i+2k},\dots$. 
This completes the induction step.
$\Cox$

\noindent{\sc Proof of Proposition}~\ref{pr:half}:
The generating polynomial for $Z$ is $\sum_k a_k z^k$ where
$a_k = (1/2) \P (X = 2k+1) + \P (X=2k) + (1/2) \P (X=2k-1)$.
Let $g(z) = (1/2) (1+z)^2 f(z)$ where $f$ is the generating polynomial
for $X$.  Then $f \in \NR$ implies $g \in \NR$.  Applying 
Theorem~\ref{th:interlace} with $g$ in place of $f$, we have
$g = g_0 + z g_1$ where $g_0 , g_1 \in \NR$.  The $z^k$
coefficient of $g_1$ is the $z^{2k+1}$ coefficient of $g$,
which we see is equal to $a_k$.  Thus $Z$ has generating 
polynomial $g_1$, which is stable.
$\Cox$

\noindent{\sc Proof of Proposition}~\ref{pr:divide}:
The generating polynomial for $Z := \lfloor X/k \rfloor$ is 
\begin{equation} \label{eq:sum}
h(y)=\sum_{i=0}^{k-1}g_i(y)
\end{equation}
where $g_0 , \ldots , g_k$ are defined from the generating polynomial
$f$ for $X$ by~\eqref{poly}.

From the proof of Theorem~\ref{th:interlace}, we see that $(-1)^mh(s_{mk})>0$ 
for each $0\leq m\leq \frac{n-k}k$ (since the smallest root is $s_{n-k}$). 
Therefore, $h$ has a root in each of the intervals
of the form $(s_{(m+1)k},s_{mk})$ for each $0\leq m\leq \frac{n-2k}k$. 
This shows that $h$ at least $\lfloor \frac nk\rfloor-1$ negative roots.
The degree of $h$ is the largest of the degrees of the $g_i$'s, 
which is the degree of $g_0$, i.e. $\lfloor \frac nk\rfloor$. 
To capture the final negative root, we observe that 
$$(-1)^{\lfloor\frac nk\rfloor}h(s_{(\lfloor\frac nk\rfloor-1)k}) < 0
   \quad\text{and }\quad(-1)^{\lfloor\frac nk\rfloor}h(s) > 0
   \text{ for large negative } s.$$
$\Cox$

We do not know the extent to which multiplication by $a$ can be 
accomplished when $a \in (0,1)$ is not a unit fraction.  The same
construction does not work.  For example, if $X$ has pgf
$$\frac 1{20}(x+1)^2(x+4),$$
then the pgf of $Y=\lfloor \frac 23 X\rfloor$ is $\frac1{20}(y^2+6y+13)$, 
which has roots $-3\pm 2\imath$.  Thus an approach analogous to the
one for unit fractions, does not work when $a = 2/3$.

\begin{question} \label{q:2/3}
Is there an $O(1)$ stable multiplication by $2/3$?
\end{question}

A solution to the following more general stable multiplication question 
would improve the hypotheses for the CLT by lowering the 
variance requirement below $M^{2/3}$.

\begin{question} \label{q:stable}
Let $\XX$ have real stable probability generating polynomial with
maximum value $M$ and let $\aa$ be a positive real vector.  
Find a stable $o(M^{1/3})$ approximation to $\aa \cdot \xx$.
\end{question}


\bibliographystyle{alpha}
\bibliography{RP}

\end{document}